\RequirePackage{fix-cm}
\documentclass[smallextended]{svjour3}       
\smartqed  
\usepackage{graphicx}
\usepackage{amssymb}

\usepackage{amsmath}
\usepackage{mathtools}
\usepackage{amsfonts}
\usepackage{mathrsfs}
\usepackage{bbm}
\usepackage{bm}

\usepackage{mathabx}

\usepackage{algorithm}
\usepackage{algorithmic}
\usepackage{booktabs}
\usepackage{multirow}
\usepackage{enumerate}
\usepackage{setspace,supertabular}
\usepackage{longtable}
\usepackage{multirow}
\usepackage{lscape}
\usepackage{url}
\usepackage{hyperref}

\usepackage[table]{xcolor}
\usepackage{array}

\newcolumntype{x}[1]{%
>{\raggedleft\hspace{0pt}}p{#1}}%

\usepackage{setspace}

\floatname{algorithm}{Procedure}

\newtheorem{assump}{Assumption}

\newtheorem{assumpB}{Assumption}

\newtheorem{assumpC}{Assumption}

\newcommand{\eps}{\varepsilon}

\newcommand{\bx}{\mathbf{x}}

\newcommand{\bt}{\mathbf{t}}

\newcommand{\zero}{\mathbf{0}}

\newcommand{\bXi}{\bm{\Xi}}

\newcommand{\bGamma}{\bm{\Gamma}}

\newcommand{\bxi}{\bm{\xi}}

\newcommand{\bLambda}{\bm{\Lambda}}

\newcommand{\E}{\mathsf{E}}

\newcommand{\Var}{\mathsf{Var}\,}
\newcommand{\Cov}{\mathsf{Cov}\,}
\newcommand{\Corr}{\mathsf{Corr}\,}
\newcommand{\prob}{\mathsf{P}}

\newcommand{\dist}{\mathscr{D}}

\newcommand{\probabi}{\mathcal{P}}

\newcommand{\Op}{\mathcal{O}_{\probabi}}

\newcommand{\ud}{\mbox{d}}

\renewcommand{\P}{\mathsf{P}} 

\newcommand{\e}{\varepsilon} 

\definecolor{seda}{gray}{0.9}

\journalname{Metrika}
\begin{document}

\title{Testing Structural Changes in Panel Data with Small Fixed Panel Size and Bootstrap\thanks{This paper was written with the support of the Czech Science Foundation project GA\v{C}R No.~P201/13/12994P. This is a~corrected version of the original paper: Pe\v{s}tov\'{a}, B. and Pe\v{s}ta, M. (2015). Testing structural changes in panel data with small fixed panel size and bootstrap. Metrika, {\bf 78}(6), 665--689, DOI 10.1007/s00184-014-0522-8, \href{http://link.springer.com/article/10.1007\%2Fs00184-014-0522-8}{link.springer.com/article/10.1007\%2Fs00184-014-0522-8}}
}
\subtitle{(Corrected version of the original paper)}

\titlerunning{Testing Structural Changes in Panel Data with Small Fixed Panel Size}        

\author{Barbora Pe\v{s}tov\'{a}  \and
        Michal Pe\v{s}ta 
}


\institute{B. Pe\v{s}tov\'{a}
           \and
           M. Pe\v{s}ta \at
              Charles University in Prague, Faculty of Mathematics and Physics, Department of Probability and Mathematical Statistics, Sokolovsk\'{a}~83, 18675 Prague~8, Czech Republic \\
              Tel.: +420-221-913-400, Fax: +420-222-323-316, \email{Michal.Pesta@mff.cuni.cz}
}

\date{Received: \today / Accepted: date}

\maketitle

\begin{abstract}
Panel data of our interest consist of a~moderate or relatively large number of panels, while the panels contain a~small number of observations. This paper establishes testing procedures to detect a~possible common change in means of the panels. To this end, we consider a~ratio type test statistic and derive its asymptotic distribution under the no change null hypothesis. Moreover, we prove the consistency of the test under the alternative. The main advantage of such an approach is that the variance of the observations neither has to be known nor estimated. On the other hand, the correlation structure is required to be calculated. To overcome this issue, a~bootstrap technique is proposed in the way of a~completely data driven approach without any tuning parameters. The validity of the bootstrap algorithm is shown. As a~by-product of the developed tests, we introduce a~common break point estimate and prove its consistency. The results are illustrated through a~simulation study. An application of the procedure to actuarial data is presented.
\end{abstract}

\section{Introduction}\label{sec:intro}
The problem of an unknown common change in means of the panels is studied here, where the panel data consist of $N$ panels and each panel contains $T$ observations over time. Various values of the change are possible for each panel at some unknown common time $\tau=1,\ldots,N$. The panels are considered to be independent, but this restriction can be weakened. In spite of that, observations within the panel are usually not independent. It is supposed that a~common unknown dependence structure is present over the panels.

Tests for change point detection in the panel data have been proposed only in case when the panel size $T$ is sufficiently large, i.e., $T$ increases over all limits from an asymptotic point of view, cf.~\cite{CHH2013} or~\cite{HH2012}. However, the change point estimation has already been studied for finite $T$ not depending on the number of panels $N$, see~\cite{Bai2010}. The remaining task is to develop testing procedures to decide whether a~common change point is present or not in the panels, while taking into account that the length $T$ of each observation regime is fixed and can be relatively small.


\subsection{Motivation}\label{subsec:mot}
Structural changes in panel data---especially \emph{common breaks in means}---are wide spread phenomena. Our primary motivation comes from non-life insurance business, where associations in many countries uniting several insurance companies collect claim amounts paid by every insurance company each year. Such a~database of cumulative claim payments can be viewed as panel data, where insurance company $i=1,\ldots,N$ provides the total claim amount $Y_{i,t}$ paid in year $t=1,\ldots,T$ into the common database. The members of the association can consequently profit from the joint database.

For the whole association it is important to know, whether a~possible change in the claim amounts occurred during the observed time horizon. Usually, the time period is relatively short, e.g., 10--15 years. To be more specific, a~widely used and very standard actuarial method for predicting future claim amounts---called chain ladder---assumes a~kind of stability of the historical claim amounts. The formal necessary and sufficient condition is derived in~\cite{PH2012}. This paper shows a~way how to test for a~possible historical instability.

\section{Panel change point model}\label{sec:model}
Let us consider the panel change point model
\begin{equation}\label{model}
Y_{i,t}=\mu_i+\delta_i\mathcal{I}\{t>\tau\}+\sigma\eps_{i,t},\quad 1\leq i\leq N,\, 1\leq t\leq T;
\end{equation}
where $\sigma>0$ is an unknown variance-scaling parameter and $T$ is fixed, not depending on $N$. The possible \emph{common change point time} is denoted by $\tau\in\{1,\ldots,T\}$. A~situation where $\tau=T$ corresponds to \emph{no change} in means of the panels. The means $\mu_i$ are panel-individual. The amount of the break in mean, which can also differ for every panel, is denoted by $\delta_i$. Furthermore, it is assumed that the sequences of panel disturbances $\{\eps_{i,t}\}_t$ are independent and within each panel the errors form a~weakly stationary sequence with a~common correlation structure. This can be formalized in the following assumption.

\begin{assump}\label{ass:A1}
\normalfont The vectors $[\eps_{i,1},\ldots,\eps_{i,T}]^{\top}$ existing on a~probability space $(\Omega,\mathcal{F},\prob)$ are $iid$ for $i=1,\ldots,N$ with $\E\eps_{i,t}=0$ and $\Var\eps_{i,t}=1$, having the autocorrelation function
$$
\rho_t=\Corr\left(\eps_{i,s},\eps_{i,s+t}\right)=\Cov\left(\eps_{i,s},\eps_{i,s+t}\right),\quad\forall s\in\{1,\ldots,T-t\},
$$
which is independent of the lag $s$, the cumulative autocorrelation function
$$
r(t)=\Var\sum_{s=1}^t \eps_{i,s}=\sum_{|s|<t}(t-|s|)\rho_s,
$$
and the shifted cumulative correlation function
$$
R(t,v)=\Cov\left(\sum_{s=1}^t\eps_{i,s},\sum_{u=t+1}^v\eps_{i,u}\right)=\sum_{s=1}^t\sum_{u=t+1}^v\rho_{u-s},\quad t<v
$$
for all $i=1,\ldots,N$ and $t,v=1,\ldots,T$.
\end{assump}

The sequence $\{\eps_{i,t}\}_{t=1}^T$ can be viewed as a~part of a~\emph{weakly stationary} process. Note that the dependent errors within each panel do not necessarily need to be linear processes. For example, GARCH processes as error sequences are allowed as well. The assumption of independent panels can indeed be relaxed, but it would make the setup much more complex. Consequently, probabilistic tools for dependent data need to be used (e.g., suitable versions of the central limit theorem). Nevertheless, assuming, that the claim amounts for different insurance companies are independent, is reasonable. Moreover, the assumption of a~common homoscedastic variance parameter $\sigma$ can be generalized by introducing weights $w_{i,t}$, which are supposed to be known. Being particular in actuarial practice, it would mean to normalize the total claim amount by the premium received, since bigger insurance companies are expected to have higher variability in total claim amounts paid.

It is required to test the \emph{null hypothesis} of no change in the means
\[
H_0:\,\tau=T
\]
against the~\emph{alternative} that at least one panel has a~change in mean
\[
H_1:\,\tau<T\quad\mbox{and}\quad\exists i\in\{1,\ldots,N\}:\,\delta_i\neq 0.
\]

\section{Test statistic and asymptotic results}\label{sec:results}
We propose a~\emph{ratio type statistic} to test $H_0$ against $H_1$, because this type of statistic does not require estimation of the nuisance parameter for the variance. Generally, this is due to the fact that the variance parameter simply cancels out from the nominator and denominator of the statistic. In spite of that, the common variance could be estimated from all the panels, of which we possess a~sufficient number. Nevertheless, we aim to construct a~valid and completely data driven testing procedure without interfering estimation and plug-in estimates instead of nuisance parameters. A~bootstrap add-on is going to serve this purpose as it is seen later on.

For surveys on ratio type test statistics, we refer to \cite{ChenTian2014}, \cite{CH1997}, \cite{HHH2008}, \cite{LZZ2008}, and~\cite{Barborka2011}. Our particular panel change point test statistic is
\[
\mathcal{R}_N(T)=\max_{t=2,\ldots,T-2}\frac{\max_{s=1,\ldots,t}\left|\sum_{i=1}^N\left[\sum_{r=1}^s\left(Y_{i,r}-\widebar{Y}_{i,t}\right)\right]\right|}{\max_{s=t,\ldots,T-1}\left|\sum_{i=1}^N\left[\sum_{r=s+1}^T\left(Y_{i,r}-\widetilde{Y}_{i,t}\right)\right]\right|},
\]
where $\widebar{Y}_{i,t}$ is the average of the first $t$ observations in panel $i$ and $\widetilde{Y}_{i,t}$ is the average of the last $T-t$ observations in panel $i$, i.e.,
\[
\widebar{Y}_{i,t}=\frac{1}{t}\sum_{s=1}^t Y_{i,s}\quad\mbox{and}\quad\widetilde{Y}_{i,t}=\frac{1}{T-t}\sum_{s=t+1}^T Y_{i,s}.
\]

An alternative way for testing the change in panel means could be a~usage of CUSUM type statistics. For example, a~maximum or minimum of a~sum (not a~ratio) of properly standardized or modified sums from our test statistic $\mathcal{R}_N(T)$. The theory, which follows, can be appropriately rewritten for such cases.

Firstly, we derive the behavior of the test statistics under the null hypothesis.

\begin{theorem}[Under null]\label{underNull}
Under hypothesis $H_0$ and Assumption~\ref{ass:A1}
\[
\mathcal{R}_N(T)\xrightarrow[N\to\infty]{\dist}\max_{t=2,\ldots,T-2}\frac{\max_{s=1,\ldots,t}\left|X_s-\frac{s}{t}X_t\right|}{\max_{s=t,\ldots,T-1}\left|Z_s-\frac{T-s}{T-t}Z_t\right|},
\]
where $Z_t:=X_T-X_t$ and $[X_1,\ldots,X_T]^{\top}$ is a~multivariate normal random vector with zero mean and covariance matrix $\bLambda=\{\lambda_{t,v}\}_{t,v=1}^{T,T}$ such that
$$
\lambda_{t,t}=r(t)\quad\mbox{and}\quad\lambda_{t,v}=r(t)+R(t,v),\,\, t<v.
$$
\end{theorem}

The limiting distribution does not depend on the variance nuisance parameter $\sigma$, but it depends on the unknown correlation structure of the panel change point model, which has to be estimated for testing purposes. The way of its estimation is shown in Subsection~\ref{subsec:covest}. Furthermore, Theorem~\ref{underNull} is just a~theoretical mid-step for the bootstrap test, where the correlation structure need not to be known. That is why the presence of unknown quantities in the asymptotic distribution is not troublesome.

Note that in case of independent observations within the panel, the correlation structure and, hence, the covariance matrix $\bLambda$ is simplified such that $r(t)=t$ and $R(t,v)=0$.

Next, we show how the test statistic behaves under the alternative.

\begin{assump}\label{alternativeDelta}
\normalfont $\lim_{N\to\infty}\frac{1}{\sqrt{N}}\left|\sum_{i=1}^N\delta_i\right|=\infty$.
\end{assump}

\begin{theorem}[Under alternative]\label{underAlternative}
If $\tau\leq T-3$, then under Assumptions~\ref{ass:A1}, \ref{alternativeDelta} and alternative $H_1$
\begin{equation}
\mathcal{R}_{N}(T)\xrightarrow[N\to\infty]{\prob}\infty.
\end{equation}
\end{theorem}

Assumption~\ref{alternativeDelta} is satisfied, for instance, if $0<\delta\leq\delta_i\,\forall i$ (a~common lower change point threshold) and $\delta\sqrt{N}\to\infty,\, N\to\infty$. Another suitable example of $\delta_i$s for the condition in Assumption~\ref{alternativeDelta}, can be $0<\delta_i=KN^{-1/2+\eta}$ for some $K>0$ and $\eta>0$. Or $\delta_i=Ci^{\alpha-1}\sqrt{N}$ may be used as well, where $\alpha\geq 0$ and $C>0$. The assumption $\tau\leq T-3$ means that there are at least three observations in the panel after the change point. It is also possible to redefine the test statistic by interchanging the nominator and the denominator of $\mathcal{R}_{N}(T)$. Afterwards, Theorem~\ref{underAlternative} for the modified test statistic would require three observations before the change point, i.e., $\tau\geq 3$.

Theorem~\ref{underAlternative} says that in presence of a~structural change in the panel means, the test statistic explodes above all bounds. Hence, the procedure is consistent and the asymptotic distribution from Theorem~\ref{underNull} can be used to construct the test.

\section{Change point estimation}
Despite the fact that the aim of the paper is to establish testing procedures for detection of a~panel mean change, it is necessary to construct a~\emph{consistent estimate} for a~possible change point. There are two reasons for that: Firstly, the estimation of the covariance matrix $\bLambda$ from Theorem~1 
requires panels as vectors with elements having common mean (i.e., without a~jump). Secondly, the bootstrap procedure, introduced later on, requires centered residuals to be resampled.

A~consistent estimate of the change point in the panel data is proposed in~\cite{Bai2010}, but under circumstances that the change occurred for sure. In our situation, we do not know whether a~change occurs or not. Therefore, we modify the estimate proposed by~\cite{Bai2010} in the following way. If the panel means change somewhere inside $\{2,\ldots,T-1\}$, let the estimate consistently select this change. If there is no change in panel means, the estimate points out the very last time point $T$ with probability going to one. In other words, the value of the change point estimate can be $T$ meaning no change. This is in contrast with~\cite{Bai2010}, where $T$ is not reachable.

Let us define the estimate of $\tau$ as
\begin{equation}\label{tauhat}
\widehat{\tau}_N:=\arg\min_{t=2,\ldots,T}\frac{1}{w(t)}\sum_{i=1}^N\sum_{s=1}^t(Y_{i,s}-\widebar{Y}_{i,t})^2,
\end{equation}
where $\{w(t)\}_{t=2}^T$ is a~sequence of weights specified later on.


\begin{assumpC}\label{ass:decreasing}
The sequence $\left\{\frac{t}{w(t)}\left(1-\frac{r(t)}{t^2}\right)\right\}_{t=2}^T$ is decreasing.
\end{assumpC}

\begin{assumpC}\label{ass:delta2E}
There exist constants $L>0$ and $N_0\in\mathbb{N}$ such that
\[
L<\sigma^2\left[\frac{t}{w(t)}\left(1-\frac{r(t)}{t^2}\right)-\frac{\tau}{w(\tau)}\left(1-\frac{r(\tau)}{\tau^2}\right)\right]+\frac{\tau(t-\tau)}{tw(t)}\frac{1}{N}\sum_{i=1}^N\delta_i^2,
\]
for each $t=\tau+1,\ldots,T$ and $N\geq N_0$.
\end{assumpC}

\begin{assumpC}\label{ass:tozero}
\normalfont $\lim_{N\to\infty}\frac{1}{N^2}\sum_{i=1}^N\delta_i^2=0$.
\end{assumpC}

\begin{assumpC}\label{ass:four}
\normalfont $\E\e_{1,t}^4<\infty,\,t\in\{1,\ldots,T\}$.
\end{assumpC}

\begin{theorem}[Change point estimate consistency]\label{thm:change}
Suppose that $\tau\neq 1$. Then under Assumptions~A1, 
\ref{ass:decreasing}, \ref{ass:delta2E}, \ref{ass:tozero}, and~\ref{ass:four}
\[
\lim_{N\to\infty}\P[\widehat{\tau}_N=\tau]=1.
\]
\end{theorem}

Assumption~\ref{ass:delta2E} assures that the values of changes have to be large enough compared to the variability of the random noise in the panels and to the strength of dependencies within the panels as well. 
Assumption~\ref{ass:tozero} is needed to control the asymptotic boundedness of the variability of $\frac{1}{w(t)}\sum_{i=1}^N\sum_{s=1}^t(Y_{i,s}-\widebar{Y}_{i,t})^2$, because a~finite $T$ cannot do that.

Assumptions~\ref{ass:delta2E} and~\ref{ass:tozero} are satisfied for $0<\delta\leq\delta_i<\Delta,\forall i$ (a~common lower and upper bound for the change amount) and suitable $\sigma$, $r(t)$, and $w(t)$. The monotonicity Assumption~\ref{ass:decreasing} in not very restrictive at all. For example in case of independent observations within the panel (i.e., $r(t)=t$) and weight function $w(t)=t^q,\,q\geq 2$, this assumption is automatically fulfilled, since sequence $\{t^{1-q}-t^{-q}\}_{t=2}^T$ is decreasing. This also gives us an~idea how to choose weights $w(t)$.

If one is interested in sensitivity of the change point estimate (i.e., what is the size of the change that can be estimated), let us consider the following model scenario: $T=10$, $\tau=5$, $\sigma=0.1$, independent observations within the panel, and $w(t)=t^2$. Then, Assumption~\ref{ass:delta2E} is satisfied if $\frac{1}{N}\sum_{i=1}^N\delta_i^2>0.029$ for all $N\geq N_0$. In case of a~common value of $\delta=\delta_i$ for all $i$, we need $\delta>\sqrt{0.029}\approx 0.170$.

Assumption~\ref{ass:delta2E} can be considered as too complicated. Therefore, one can replace it by the following simpler, but more restrictive assumption.
\begin{assumpC}\label{ass:infty}
\[
\lim_{N\to\infty}\frac{1}{N}\sum_{i=1}^N\delta_i^2=\infty.
\]
\end{assumpC}

On one hand, this assumption might be considered as too strong, because a~common fixed (not depending on $N$) value of $\delta=\delta_i$ for all $i$ does not fulfill Assumption~\ref{ass:infty}. On the other hand, Assumption~\ref{ass:infty} is satisfied when $\delta_j^2/N\to\infty$ as $N\to\infty$ for some $j\in\mathbb{N}$ and $\delta_i=0$ for all $i\neq j$. This stands for a~situation when all the panels do not change in mean except one panel having a~sufficiently large change in mean with respect to the number of panels.

Various competing consistent estimates of a~possible change point can be suggested, e.g., the maximizer of $\sum_{i=1}^N\left[\sum_{s=1}^t(Y_{i,s}-\widebar{Y}_{i,T})\right]^2$. To show the consistency, one needs to postulate different assumptions on the cumulative autocorrelation function and shifted cumulative correlation function compared to Theorem~\ref{thm:change} and this may be rather complex.

\subsection{Estimation of the correlation structure}\label{subsec:covest}
Since the panels are considered to be independent and the number of panels may be sufficiently large, one can estimate the correlation structure of the errors $[\eps_{1,1},\ldots,\eps_{1,T}]^{\top}$ empirically. We base the errors' estimates on \emph{residuals}
\begin{equation}\label{widehate}
\widehat{e}_{i,t}:=\left\{
\begin{array}{ll}
Y_{i,t}-\widebar{Y}_{i,\widehat{\tau}_N},& t\leq\widehat{\tau}_N,\\
Y_{i,t}-\widetilde{Y}_{i,\widehat{\tau}_N},& t>\widehat{\tau}_N.
\end{array}
\right.
\end{equation}

Then, the empirical version of the autocorrelation function is
\[
\widehat{\rho}_t:=\frac{1}{\widehat{\sigma}^2 NT}\sum_{i=1}^N\sum_{s=1}^{T-t}\widehat{e}_{i,s}\widehat{e}_{i,s+t}.
\]
Consequently, the kernel estimation of the cumulative autocorrelation function and shifted cumulative correlation function is adopted in lines with~\cite{Andrews1991}:
\begin{align*}
\widehat{r}(t)&=\sum_{|s|<t}(t-|s|)\kappa\left(\frac{s}{h}\right)\widehat{\rho}_s,\\
\widehat{R}(t,v)&=\sum_{s=1}^t\sum_{u=t+1}^v\kappa\left(\frac{u-s}{h}\right)\widehat{\rho}_{u-s},\quad t<v;
\end{align*}
where $h>0$ stands for the window size and $\kappa$ belongs to a~class of kernels given by
\begin{multline*}
\Big\{\kappa(\cdot):\,\mathbbm{R}\to[-1,1]\,\big|\,\kappa(0)=1,\,\kappa(x)=\kappa(-x),\,\forall x,\,\int_{-\infty}^{+\infty}\kappa^2(x)\ud x<\infty,\Big.\\
\Big.\kappa(\cdot)\,\mbox{ is continuos at $0$ and at all but a finite number of other points}\Big\}.
\end{multline*}

Since the variance parameter $\sigma$ is not present in the limiting distribution of Theorem~\ref{underNull}, it neither has to be estimated nor known. Nevertheless, one can use $\widehat{\sigma}^2:=\frac{1}{NT}\sum_{i=1}^{N}\sum_{s=1}^{T}\widehat{e}_{i,s}^2$.

\section{Bootstrap and hypothesis testing}\label{sec:boot}
A~wide range of literature has been published on bootstrapping in the change point problem, e.g., \cite{HK2012} or~\cite{HKPS2008}. We build up the bootstrap test on the resampling with replacement of row vectors $\{[\widehat{e}_{i,1},\ldots,\widehat{e}_{i,T}]\}_{i=1,\ldots,N}$ corresponding to the panels. This provides bootstrapped row vectors $\{[\widehat{e}_{i,1}^*,\ldots,\widehat{e}_{i,T}^*]\}_{i=1,\ldots,N}$. Then, the bootstrapped residuals $\widehat{e}_{i,t}^*$ are \emph{centered} by their conditional expectation $\frac{1}{N}\sum_{i=1}^N\widehat{e}_{i,t}$ yielding
\[
\widehat{Y}_{i,t}^*:=\widehat{e}_{i,t}^*-\frac{1}{N}\sum_{i=1}^N\widehat{e}_{i,t}.
\]
The bootstrap test statistic is just a~modification of the original statistic $\mathcal{R}_N(T)$, where the original observations $Y_{i,t}$ are replaced by their bootstrap counterparts $\widehat{Y}_{i,t}^*$:
\[
\mathcal{R}_{N}^*(T)=\max_{t=2,\ldots,T-2}\frac{\max_{s=1,\ldots,t}\left|\sum_{i=1}^N\left[\sum_{r=1}^s\left(\widehat{Y}_{i,r}^*-\widebar{\widehat{Y}}_{i,t}^*\right)\right]\right|}{\max_{s=t,\ldots,T-1}\left|\sum_{i=1}^N\left[\sum_{r=s+1}^T\left(\widehat{Y}_{i,r}^*-\widetilde{\widehat{Y}}_{i,t}^*\right)\right]\right|},
\]
such that
\[
\widebar{\widehat{Y}}_{i,t}^*=\frac{1}{t}\sum_{s=1}^t \widehat{Y}_{i,s}^*\quad\mbox{and}\quad\widetilde{\widehat{Y}}_{i,t}^*=\frac{1}{T-t}\sum_{s=t+1}^T \widehat{Y}_{i,s}^*.
\]

An~\emph{algorithm} for the bootstrap is illustratively shown in Procedure~\ref{alg:nonpar} and its validity will be proved in Theorem~\ref{bootJust}.
\begin{algorithm}[!htb]
\caption{Bootstrapping test statistic $\mathcal{R}_N(T)$.}
\label{alg:nonpar}
\begin{algorithmic}[1]
\REQUIRE Panel data consisting of $N$ panels with length $T$, i.e., $N$ row vectors of observations $[Y_{i,1},\ldots,Y_{i,T}]$.
\ENSURE Bootstrap distribution of $\mathcal{R}_N(T)$, i.e., the empirical distribution where probability mass $1/B$ concentrates at each of ${}_{(1)}\mathcal{R}_{N}^*(T),\ldots,{}_{(B)}\mathcal{R}_{N}^*(T)$.
\STATE estimate the change point by calculating $\widehat{\tau}_N$
\STATE compute residuals $\widehat{e}_{i,t}$
\FOR[repeat in order to obtain the empirical distribution]{$b=1$ to $B$}
\STATE $\{[\widehat{e}_{i,1}^*,\ldots,\widehat{e}_{i,T}^*]\}_{i=1}^N$ resampled with replacement from original rows $\{[\widehat{e}_{i,1},\ldots,\widehat{e}_{i,T}]\}_{i=1}^N$
\STATE calculate bootstrap panel data $\widehat{Y}_{i,t}^*$
\STATE compute bootstrap test statistics ${}_{(b)}\mathcal{R}_{N}^*(T)$
\ENDFOR
\end{algorithmic}
\end{algorithm}

\subsection{Validity of the resampling procedure}\label{subsec:val}
The idea behind bootstrapping is to \emph{mimic the original distribution} of the test statistic in some sense with the distribution of the bootstrap test statistic, conditionally on the original data denoted by $\mathbbm{Y}\equiv\{Y_{i,t}\}_{i,t=1}^{N,T}$.

First of all, two simple and just technical assumptions are needed.

\begin{assumpB}\label{lagged}
\normalfont $\{\eps_{i,t}\}_t$ possesses the lagged cumulative correlation function
$$
S(t,v,d)=\Cov\left(\sum_{s=1}^t\eps_{i,s},\sum_{u=t+d}^v\eps_{i,u}\right)=\sum_{s=1}^t\sum_{u=t+d}^v\rho_{u-s},\quad\forall i\in\mathbbm{N}.
$$
\end{assumpB}

\begin{assumpB}\label{withprob}
\normalfont $\lim_{N\to\infty}\prob[\widehat{\tau}_N=\tau]=1$.
\end{assumpB}

Assumption~\ref{lagged} is not really an~assumption, actually it is only a~notation. Notice that $S(t,v,1)\equiv R(t,v)$. Assumption~\ref{withprob} is satisfied for our estimate proposed in~\eqref{tauhat}, if the assumptions of Theorem~\ref{underAlternative} hold. Assumption~\ref{withprob} is postulated in a~rather broader sense, because we want to allow any other consistent estimate of $\tau$ to be used instead.

Realize that it is not known, whether the common panel means' change occurred or not. In other words, one does not know \emph{whether the data come from the null or the alternative} hypothesis. Therefore, the following theorem holds under $H_0$ as well as $H_1$.

\begin{theorem}[Bootstrap justification]\label{bootJust}
Under Assumptions~\ref{ass:A1}, \ref{lagged}, \ref{withprob}, and~\ref{ass:four}
\begin{equation*}
\mathcal{R}_{N}^*(T)|\mathbbm{Y}\xrightarrow[N\to\infty]{\dist}\max_{t=2,\ldots,T-2}\frac{\max_{s=1,\ldots,t}\left|\mathcal{X}_s-\frac{s}{t}\mathcal{X}_t\right|}{\max_{s=t,\ldots,T-1}\left|\mathcal{Z}_s-\frac{T-s}{T-t}\mathcal{Z}_t\right|}\quad\mbox{in probability }\prob,
\end{equation*}
where $\mathcal{Z}_t:=\mathcal{X}_T-\mathcal{X}_t$ and $[\mathcal{X}_1,\ldots,\mathcal{X}_T]^{\top}$ is a~multivariate normal random vector with zero mean and covariance matrix $\bGamma=\left\{\gamma_{t,v}(\tau)\right\}_{t,v=1}^{T,T}$ such that
\begin{equation*}
\gamma_{t,t}(\tau)=\left\{\begin{array}{l}
r(t)+\frac{t^2}{\tau^2}r(\tau)-\frac{2t}{\tau}[r(t)+R(t,\tau)],\\
\quad t<\tau;\\
0,\quad t=\tau;\\
r(t-\tau)+\frac{(t-\tau)^2}{(T-\tau)^2} r(T-\tau)-\frac{2(t-\tau)}{T-\tau}\left[r(t-\tau)+R(t-\tau,T-\tau)\right],\\
\quad t>\tau;
\end{array}
\right.
\end{equation*}
and
\begin{equation*}
\gamma_{t,v}(\tau)=\left\{\begin{array}{l}
0,\quad t=\tau\mbox{ or }v=\tau,\\
r(t)+R(t,v)+\frac{tv}{\tau^2}r(\tau)-\frac{v}{\tau}[r(t)+R(t,\tau)]\\
\quad -\frac{t}{\tau}[r(v)+R(v,\tau)],\quad t<v<\tau;\\
S(t,v,\tau+1-t)+\frac{t(v-\tau)}{\tau(T-\tau)}R(\tau,T)\\
\quad -\frac{v-\tau}{T-\tau}S(t,T,\tau+1-t)-\frac{t}{\tau}R(\tau,v),\quad t<\tau<v;\\
r(t-\tau)+R(t-\tau,v-\tau)+\frac{(t-\tau)(v-\tau)}{(T-\tau)^2}r(T-\tau)\\
\quad -\frac{v-\tau}{T-\tau}[r(t-\tau)+R(t-\tau,T-\tau)]\\
\quad -\frac{t-\tau}{T-\tau}[r(v-\tau)+R(v-\tau,T-\tau)],\quad \tau<t<v.
\end{array}
\right.
\end{equation*}
\end{theorem}

The validity of the bootstrap test is assured by Theorem~\ref{bootJust}. Indeed, the conditional asymptotic distribution of the bootstrap test statistic is a~functional of a~multivariate normal distribution under the null as well as under the alternative. It does not converge to infinity (in probability) under the alternative. That is why it can be used for correctly rejecting the null in favor of the alternative, having sufficiently large $N$. Moreover, the following theorem states that the conditional distribution of the bootstrap test statistic and the unconditional distribution of the original test statistic \emph{coincide}. And that is the reason why the bootstrap test should approximately keep the same level as the original test based on the asymptotics from Theorem~\ref{underNull}.

\begin{theorem}[Bootstrap test consistency]\label{bootConsis}
Under Assumptions~\ref{ass:A1}, \ref{withprob}, \ref{ass:four} and hypothesis $H_0$, the asymptotic distribution of $\mathcal{R}_{N}(T)$ from Theorem~\ref{underNull} and the asymptotic distribution of $\mathcal{R}_{N}^*(T)|\mathbbm{Y}$ from Theorem~\ref{bootJust} coincide.
\end{theorem}

Now, the simulated (empirical) distribution of the bootstrap test statistic can be used to calculate the bootstrap critical value, which will be compared to the value of the original test statistic in order to reject the null or not.

Finally, note that one cannot think about any local alternative in this setup, because $\tau$ has a~discrete and finite support.

\section{Simulations}\label{sec:simul}
A~simulation experiment was performed to study the \emph{finite sample} properties of the asymptotic and bootstrap test statistics for a~common change in panel means. In particular, the interest lies in the empirical \emph{sizes} of the proposed tests under the null hypothesis and in the empirical \emph{rejection} rate (power) under the alternative. Random samples of panel data ($5000$ each time) are generated from the panel change point model~\eqref{model}. The panel size is set to $T=10$ and $T=25$ in order to demonstrate the performance of the testing approaches in case of small and intermediate panel length. The number of panels considered is $N=50$ and $N=200$.

The correlation structure within each panel is modeled via random vectors generated from iid, AR(1), and GARCH(1,1) sequences. The considered AR(1) process has coefficient $\phi=0.3$. In case of GARCH(1,1) process, we use coefficients $\alpha_0=1$, $\alpha_1=0.1$, and $\beta_1=0.2$, which according to \cite[Example~1]{Lindner2009} gives a~strictly stationary process. In all three sequences, the innovations are obtained as iid random variables from a~standard normal $\mathsf{N}(0,1)$ or Student $t_5$ distribution. Simulation scenarios are produced as all possible combinations of the above mentioned settings.

When using the asymptotic distribution from Theorem~\ref{underNull}, the covariance matrix is estimated as proposed in Subsection~\ref{subsec:covest} using the Parzen kernel
\[
\kappa_{P}(x)=\left\{\begin{array}{ll}
1-6x^2+6|x|^3, & 0\leq|x|\leq 1/2;\\
2(1-|x|)^3, & 1/2\leq|x|\leq 1;\\
0, & \mbox{otherwise}.
\end{array}
\right.
\]
Several values of the smoothing window width $h$ are tried from the interval $[2,5]$ and all of them work fine providing comparable results. To simulate the asymptotic distribution of the test statistics, $2000$ multivariate random vectors are generated using the pre-estimated covariance matrix.

The bootstrap approach does not need to estimate the covariance structure. The number of bootstrap replications used is $2000$. To access the theoretical results under $H_0$ numerically, Table~\ref{tab:H0} provides the empirical specificity (one minus size) of the tests for both the asymptotic and bootstrap version of the panel change point test, where the significance level is $\alpha=5\%$.
\begin{table}[!ht]
\caption{Empirical specificity ($1-$size) of the test under $H_0$ using the asymptotic and the \colorbox{seda}{bootstrap} critical values, considering a~significance level of $5\%$, weight function $w(t)=t^2$, and smoothing window width $h=2$.}
\label{tab:H0}
\begin{tabular}{cccx{.65cm}x{.65cm}x{.65cm}x{.65cm}x{.65cm}x{.65cm}}
\toprule
$T$ & $N$ & innovations & \multicolumn{2}{c}{IID} & \multicolumn{2}{c}{AR(1)} & \multicolumn{2}{c}{GARCH(1,1)} \tabularnewline
\midrule
$10$ & $50$ & $\mathsf{N}(0,1)$ & $.942$ & \cellcolor[gray]{0.9}{$.959$} & $.932$ & \cellcolor[gray]{0.9}{$.962$} & $.952$ & \cellcolor[gray]{0.9}{$.968$} \tabularnewline
 & & $t_5$ & $.950$ & \cellcolor[gray]{0.9}{$.967$} & $.933$ & \cellcolor[gray]{0.9}{$.962$} & $.947$ & \cellcolor[gray]{0.9}{$.966$} \tabularnewline
 & $200$ & $\mathsf{N}(0,1)$ & $.950$ & \cellcolor[gray]{0.9}{$.964$} & $.938$ & \cellcolor[gray]{0.9}{$.968$} & $.950$ & \cellcolor[gray]{0.9}{$.968$} \tabularnewline
 & & $t_5$ & $.950$ & \cellcolor[gray]{0.9}{$.964$} & $.934$ & \cellcolor[gray]{0.9}{$.964$} & $.941$ & \cellcolor[gray]{0.9}{$.963$} \tabularnewline
\cmidrule(){1-9}
$25$ & $50$ & $\mathsf{N}(0,1)$ & $.945$ & \cellcolor[gray]{0.9}{$.961$} & $.933$ & \cellcolor[gray]{0.9}{$.965$} & $.947$ & \cellcolor[gray]{0.9}{$.963$} \tabularnewline
 & & $t_5$ & $.949$ & \cellcolor[gray]{0.9}{$.964$} & $.929$ & \cellcolor[gray]{0.9}{$.964$} & $.947$ & \cellcolor[gray]{0.9}{$.963$} \tabularnewline
 & $200$ & $\mathsf{N}(0,1)$ & $.951$ & \cellcolor[gray]{0.9}{$.962$} & $.928$ & \cellcolor[gray]{0.9}{$.964$} & $.953$ & \cellcolor[gray]{0.9}{$.968$} \tabularnewline
 & & $t_5$ & $.954$ & \cellcolor[gray]{0.9}{$.965$} & $.931$ & \cellcolor[gray]{0.9}{$.966$} & $.953$ & \cellcolor[gray]{0.9}{$.967$} \tabularnewline
\bottomrule
\end{tabular}
\end{table}

It may be seen that both approaches (using asymptotic and bootstrap distribution) are close to the theoretical value of specificity $.95$. As expected, the best results are achieved in case of independence within the panel, because there is no information overlap between two consecutive observations. The precision of not rejecting the null is increasing as the number of panels is getting higher and the panel is getting longer as well.

The performance of both testing procedures under $H_1$ in terms of the empirical rejection rates is shown in Table~\ref{tab:H1}, where the change point is set to $\tau=\lfloor T/2 \rfloor$ and the change sizes $\delta_i$ are independently uniform on $[1,3]$ in $33\%$, $66\%$ or in all panels.
\begin{table}[!ht]
\caption{Empirical sensitivity (power) of the test under $H_1$ using the asymptotic and the \colorbox{seda}{bootstrap} critical values, considering a~significance level of $5\%$, weight function $w(t)=t^2$, and smoothing window width $h=2$.}
\label{tab:H1}
\begin{tabular}{ccccx{.65cm}x{.65cm}x{.65cm}x{.65cm}x{.65cm}x{.65cm}}
\toprule
$H_1$ & $T$ & $N$ & innovations & \multicolumn{2}{c}{IID} & \multicolumn{2}{c}{AR(1)} & \multicolumn{2}{c}{GARCH(1,1)} \tabularnewline
\midrule
$33\%$ & $10$ & $50$ & $\mathsf{N}(0,1)$ & $.23$ & \cellcolor[gray]{0.9}{$.06$} & $.26$ & \cellcolor[gray]{0.9}{$.07$} & $.19$ & \cellcolor[gray]{0.9}{$.05$} \tabularnewline
& & & $t_5$ & $.18$ & \cellcolor[gray]{0.9}{$.05$} & $.20$ & \cellcolor[gray]{0.9}{$.06$} & $.20$ & \cellcolor[gray]{0.9}{$.05$} \tabularnewline
& & $200$ & $\mathsf{N}(0,1)$ & $.45$ & \cellcolor[gray]{0.9}{$.05$} & $.48$ & \cellcolor[gray]{0.9}{$.05$} & $.39$ & \cellcolor[gray]{0.9}{$.05$} \tabularnewline
& & & $t_5$ & $.36$ & \cellcolor[gray]{0.9}{$.05$} & $.39$ & \cellcolor[gray]{0.9}{$.05$} & $.39$ & \cellcolor[gray]{0.9}{$.05$} \tabularnewline
\cmidrule(l){2-10}
& $25$ & $50$ & $\mathsf{N}(0,1)$ & $.38$ & \cellcolor[gray]{0.9}{$.05$} & $.39$ & \cellcolor[gray]{0.9}{$.05$} & $.31$ & \cellcolor[gray]{0.9}{$.05$} \tabularnewline
& & & $t_5$ & $.30$ & \cellcolor[gray]{0.9}{$.05$} & $.30$ & \cellcolor[gray]{0.9}{$.05$} & $.31$ & \cellcolor[gray]{0.9}{$.06$} \tabularnewline
& & $200$ & $\mathsf{N}(0,1)$ & $.68$ & \cellcolor[gray]{0.9}{$.05$} & $.70$ & \cellcolor[gray]{0.9}{$.05$} & $.58$ & \cellcolor[gray]{0.9}{$.05$} \tabularnewline
& & & $t_5$ & $.56$ & \cellcolor[gray]{0.9}{$.05$} & $.57$ & \cellcolor[gray]{0.9}{$.05$} & $.59$ & \cellcolor[gray]{0.9}{$.05$} \tabularnewline
\cmidrule(){1-10}
$66\%$ & $10$ & $50$ & $\mathsf{N}(0,1)$ & $.45$ & \cellcolor[gray]{0.9}{$.39$} & $.49$ & \cellcolor[gray]{0.9}{$.46$} & $.38$ & \cellcolor[gray]{0.9}{$.10$} \tabularnewline
& & & $t_5$ & $.36$ & \cellcolor[gray]{0.9}{$.10$} & $.37$ & \cellcolor[gray]{0.9}{$.14$} & $.39$ & \cellcolor[gray]{0.9}{$.15$} \tabularnewline
& & $200$ & $\mathsf{N}(0,1)$ & $.77$ & \cellcolor[gray]{0.9}{$.59$} & $.81$ & \cellcolor[gray]{0.9}{$.93$} & $.68$ & \cellcolor[gray]{0.9}{$.05$} \tabularnewline
& & & $t_5$ & $.64$ & \cellcolor[gray]{0.9}{$.05$} & $.69$ & \cellcolor[gray]{0.9}{$.12$} & $.69$ & \cellcolor[gray]{0.9}{$.10$} \tabularnewline
\cmidrule(l){2-10}
& $25$ & $50$ & $\mathsf{N}(0,1)$ & $.69$ & \cellcolor[gray]{0.9}{$.08$} & $.70$ & \cellcolor[gray]{0.9}{$.11$} & $.58$ & \cellcolor[gray]{0.9}{$.05$} \tabularnewline
& & & $t_5$ & $.56$ & \cellcolor[gray]{0.9}{$.05$} & $.57$ & \cellcolor[gray]{0.9}{$.06$} & $.59$ & \cellcolor[gray]{0.9}{$.06$} \tabularnewline
& & $200$ & $\mathsf{N}(0,1)$ & $.95$ & \cellcolor[gray]{0.9}{$.06$} & $.96$ & \cellcolor[gray]{0.9}{$.05$} & $.91$ & \cellcolor[gray]{0.9}{$.05$} \tabularnewline
& & & $t_5$ & $.87$ & \cellcolor[gray]{0.9}{$.05$} & $.89$ & \cellcolor[gray]{0.9}{$.05$} & $.91$ & \cellcolor[gray]{0.9}{$.05$} \tabularnewline
\cmidrule(){1-10}
$100\%$ & $10$ & $50$ & $\mathsf{N}(0,1)$ & $.64$ & \cellcolor[gray]{0.9}{$.92$} & $.67$ & \cellcolor[gray]{0.9}{$.84$} & $.56$ & \cellcolor[gray]{0.9}{$.58$} \tabularnewline
& & & $t_5$ & $.52$ & \cellcolor[gray]{0.9}{$.37$} & $.55$ & \cellcolor[gray]{0.9}{$.45$} & $.55$ & \cellcolor[gray]{0.9}{$.55$} \tabularnewline
& & $200$ & $\mathsf{N}(0,1)$ & $.93$ & \cellcolor[gray]{0.9}{$1.00$} & $.95$ & \cellcolor[gray]{0.9}{$1.00$} & $.87$ & \cellcolor[gray]{0.9}{$.77$} \tabularnewline
& & & $t_5$ & $.85$ & \cellcolor[gray]{0.9}{$.36$} & $.87$ & \cellcolor[gray]{0.9}{$.74$} & $.87$ & \cellcolor[gray]{0.9}{$.72$} \tabularnewline
\cmidrule(l){2-10}
& $25$ & $50$ & $\mathsf{N}(0,1)$ & $.87$ & \cellcolor[gray]{0.9}{$.84$} & $.88$ & \cellcolor[gray]{0.9}{$.85$} & $.79$ & \cellcolor[gray]{0.9}{$.11$} \tabularnewline
& & & $t_5$ & $.76$ & \cellcolor[gray]{0.9}{$.08$} & $.77$ & \cellcolor[gray]{0.9}{$.11$} & $.79$ & \cellcolor[gray]{0.9}{$.20$} \tabularnewline
& & $200$ & $\mathsf{N}(0,1)$ & $1.00$ & \cellcolor[gray]{0.9}{$.97$} & $1.00$ & \cellcolor[gray]{0.9}{$.97$} & $.99$ & \cellcolor[gray]{0.9}{$.05$} \tabularnewline
& & & $t_5$ & $.98$ & \cellcolor[gray]{0.9}{$.05$} & $.98$ & \cellcolor[gray]{0.9}{$.05$} & $.99$ & \cellcolor[gray]{0.9}{$.07$} \tabularnewline
\bottomrule
\end{tabular}
\end{table}

One can conclude that the power of both tests increases as the panel size and the number of panels increase, which is straightforward and expected. It should be noticed that numerical instability issues may appear for larger $T$, when generating from a~$T$-variate normal distribution. Moreover, higher power is obtained when a~larger portion of panels is subject to have a~change in mean. The test power drops when switching from independent observations within the panel to dependent ones. Innovations with heavier tails (i.e., $t_5$) yield smaller power than innovations with lighter tails. Generally, the bootstrap outperforms the classical asymptotics in all scenarios.

Let us mention that for finite sections of processes with a~stronger dependence structure than taken into account in the simulation scenarios, Assumption~\ref{ass:delta2E} does not have to be fulfilled. The dependency under the considered variability can be too strong compared to the change size. Then, it would be rather difficult to detect possible changes.

Finally, an~early change point is discussed very briefly. We stay with standard normal innovations, iid observations within the panel, the size of changes $\delta_i$ being independently uniform on $[1,3]$ in all panels, and the change point is $\tau=3$ in case of $T=10$ and $\tau=5$ for $T=25$. The empirical sensitivities of both tests for small values of $\tau$ are shown in Table~\ref{tab:tau}.
\begin{table}[!ht]
\caption{Empirical sensitivity of the test for small values of $\tau$ under $H_1$ using the asymptotic and the \colorbox{seda}{bootstrap} critical values, considering a~significance level of $5\%$, weight function $w(t)=t^2$, and smoothing window width $h=2$.}
\label{tab:tau}
\begin{tabular}{cccx{.65cm}x{.65cm}}
\toprule
$T$ & $N$ & $\tau$ & \multicolumn{2}{c}{$H_1$, iid, $\mathsf{N}(0,1)$} \tabularnewline
\midrule
$10$ & $50$ & $3$  & $.56$ & \cellcolor[gray]{0.9}{$.08$} \tabularnewline
 & $200$ & $3$  & $.87$ & \cellcolor[gray]{0.9}{$.05$} \tabularnewline
\cmidrule(){1-5}
$25$ & $50$ & $5$ & $.63$ & \cellcolor[gray]{0.9}{$.05$} \tabularnewline
 & $200$ & $5$  & $.92$ & \cellcolor[gray]{0.9}{$.05$} \tabularnewline
\bottomrule
\end{tabular}
\end{table}

When the change point is not in the middle of the panel, the power of the test generally falls down. The source of such decrease is that the left or right part of the panel possesses less observations with constant mean, which leads to a~decrease of precision in the correlation estimation in case of the asymptotic test and in the change point estimation in case of the bootstrap test. Nevertheless, the bootstrap test again outperforms the asymptotic version and, moreover, provides solid results even for early or late change points (the late change points are not numerically demonstrated here).

\section{Real data analysis}\label{sec:data}
\begin{figure}[!ht]
  \includegraphics[width=0.99\textwidth]{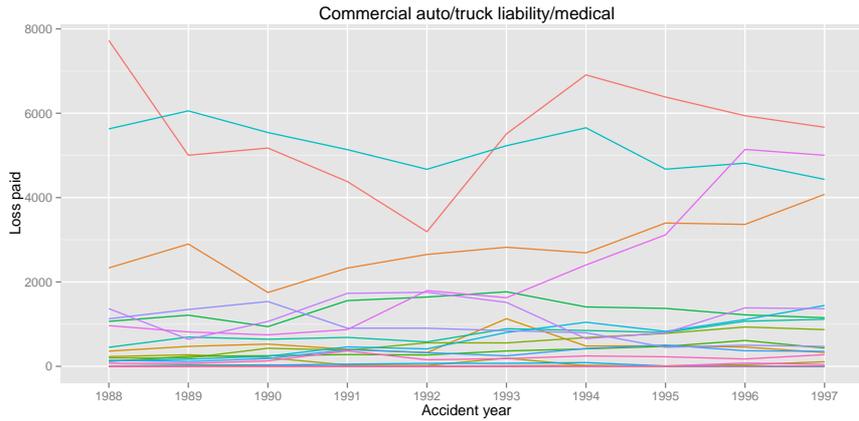}
\caption{Development of yearly total claim amounts for 20 selected insurance companies}
\label{fig:comauto20}       
\end{figure}

As mentioned in the introduction, our primary motivation for testing the panel mean change comes from the \emph{insurance business}. The data set is provided by the National Association of Insurance Commissioners (NAIC) database, see~\cite{NAIC2011}. We concentrate on the `Commercial auto/truck liability/medical' insurance line of business. The data collect records from $N=157$ insurance companies (one extreme insurance company was omitted from the analysis). Each insurance company provides $T=10$ yearly total claim amounts starting from year 1988 up to year 1997. Figure~\ref{fig:comauto20} graphically shows series of claim amounts for 20 selected insurance companies (a~plot with all 157 panels would be cluttered).

The data are considered as panel data in the way that each insurance company corresponds to one panel, which is formed by the company's yearly total claim amounts. The length of the panel is quite short. This is very typical in insurance business, because considering longer panels may invoke incomparability between the early claim amounts and the late ones due to changing market or policies' conditions over time.

We want to test whether or not a~change in the claim amounts occurred in a~common year, assuming that the claim amounts are approximately constant in the years before and after the possible change for every insurance company. Our ratio type test statistic gives $\mathcal{R}_{157}(10)=39.9$. The asymptotic critical value is $52.4$ and the bootstrap critical value equals $203.1$. These values mean that we do not reject the hypothesis of no change in panel means in both cases. The striking difference between the two critical values may come from the inefficient correlation structure estimation (since $T=10$ is quite short) or from violation of the model assumptions.

That is why we also try to take the logarithms of claim amounts and to consider log amounts as the panel data observations. Nevertheless, we do not reject the hypothesis of no change in the panel means (i.e., means of log amounts) again. Additionally to that, one can consider normalizing the claim amounts by the premium received by company $i$ in year $t$. That is thinking of panel data $Y_{i,t}/p_{i,t}$, where $p_{i,t}$ is the mentioned premium. This may yield a~stabilization of series' variability, which corresponds to the assumption of a~common variance. In spite of that, we again do not reject the null (neither by the asymptotic test, nor by the bootstrap one). For the sake of completeness, we may reveal that our estimate of the panel change point provides value $\widehat{\tau}_N=10$ meaning no change in panels.

\section{Conclusions}\label{sec:concl}
In this paper, we consider the change point problem in panel data with fixed panel size. Occurrence of common breaks in panel means is tested. We introduce a~ratio type test statistic and derive its asymptotic properties. Under the null hypothesis of no change, the test statistic weakly converges to a~functional of the multivariate normal random vector with zero mean and covariance structure depending on the intra-panel covariances. As shown in the paper, these covariances can be estimated and, consequently, used for testing whether a~change in means occurred or not. This is indeed feasible, because the test statistic under the alternative converges to infinity in probability.

The secondary aim of the paper lies in proposing a~consistent change point estimate, which is straightforwardly used for bootstrapping the test statistic. We establish the asymptotic behavior of the bootstrap version of the test statistic, regardless of the fact whether the data come from the null or the alternative hypothesis. Moreover, the asymptotic distribution of the bootstrap test statistic coincides with the original test statistic's limiting distribution. This provides justification for the bootstrap method. One of the main goals is to obtain a~completely data driven testing approach whether the means remain the same during the observation period or not. The ratio type test statistic allows us to omit a~variance estimation and the bootstrap technique overcomes estimation of the correlation structure. Hence, neither nuisance nor smoothing parameters are present in the whole testing process, which makes it very simple for practical use. Furthermore, the whole stochastic theory behind requires relatively simple assumptions, which are not too restrictive.

A~simulation study illustrates that even for small panel size, both presented approaches---based on traditional asymptotics and on bootstrapping---work fine. One may judge that both methods keep the significance level under the null, while various simulation scenarios are considered. Besides that, the power of the test is slightly higher in case of the bootstrap. Finally, the proposed methods are applied to insurance data, for which the change point analysis in panel data provides an~appealing approach.

\subsection{Discussion}\label{subsec:disc}
First of all, it has to be noted that the non-ratio CUSUM type test statistic can be used instead of the ratio type, but this requires to estimate the variance of the observations. The statements of theorems and proofs would become even less complicated. Omitting the usage of the bootstrap test statistic can especially be unreliable in short panels from a~computational point of view. This is due to the fact that the bootstrap overcomes the issue of estimating the correlation structure.

Furthermore, our setup can be modified by considering large panel size, i.e., $T\to\infty$. Consequently, the whole theory leads to convergences to functionals of Gaussian processes with a~covariance structure derived in a~very similar fashion as for fixed $T$. However, our motivation is to develop techniques for fixed and relatively small panel size.

Dependent panels may be taken into account and the presented work might be generalized for some kind of asymptotic independence over the panels or prescribed dependence among the panels. Nevertheless, our incentive is determined by a~problem from non-life insurance, where the association of insurance companies consists of a~relatively high number of insurance companies. Thus, the portfolio of yearly claims is so diversified, that the panels corresponding to insurance companies' yearly claims may be viewed as independent and neither natural ordering nor clustering has to be assumed.


\appendix



\setcounter{theorem}{0}
\renewcommand{\thetheorem}{\Alph{section}.\arabic{theorem}}

\setcounter{definition}{0}
\renewcommand{\thedefinition}{\Alph{section}.\arabic{definition}}

\section{Supporting Theorems}\label{sec:support}

Suppose that $\{\bxi_n\}_{n=1}^{\infty}$ is a~sequence of random variables/vectors existing on a~probability space $(\Omega,\mathcal{F},\prob)$. A~\emph{bootstrap version} of $\bxi\equiv[\bxi_1,\ldots,\bxi_n]^{\top}$ is its (randomly) resampled sequence with replacement---denoted by $\bxi^*\equiv[\bxi_1^*,\ldots,\bxi_n^*]^{\top}$---with the same length, where for each $i\in\{1,\ldots,n\}$ it holds that $\prob_{\bxi}^*[\bxi_i^*=\bxi_j]\equiv\prob[\bxi_i^*=\bxi_j|\bxi]=1/n,\,j=1,\ldots,n$. In the sequel, $\prob_{\bxi}^*$ denotes the conditional probability given ${\bxi}$. So, $\bxi_i^*$ has a~discrete uniform distribution on $\{\bxi_1,\ldots,\bxi_n\}$ for every $i=1,\ldots,n$. The conditional expectation and variance given ${\bxi}$ are denoted by $\E_{\prob_{\bxi}^*}$ and $\Var_{\prob_{\bxi}^*}$.

If a~statistic has an~approximate normal distribution, one may be interested in the asymptotic comparison of the bootstrap distribution with the original one. A~tool for assessing such an~approximate closeness can be a~\emph{bootstrap central limit theorem} for triangular arrays.


\begin{theorem}[Bootstrap CLT for triangular arrays]\label{thm:CLTbootTRIA}
Let $\{\xi_{n,k_n}\}_{n=1}^{\infty}$ be a~triangular array of zero mean random variables on the same probability space such that the elements of the vector $[\xi_{n,1},\ldots,\xi_{n,k_n}]^{\top}$ are iid for every $n\in\mathbbm{N}$ satisfying
\begin{equation}\label{CLTbootIID1}
\sup_{n\in\mathbbm{N}}\E_{\prob}\xi_{n,1}^4<\infty
\end{equation}
and $k_n\to\infty$ as $n\to\infty$. Suppose that $\bxi^*\equiv[\xi_{n,1}^*,\ldots,\xi_{n,k_n}^*]^{\top}$ is the bootstrapped version of $\bxi\equiv[\xi_{n,1},\ldots,\xi_{n,k_n}]^{\top}$ and denote
\[
\bar{\xi}_n:=k_n^{-1}\sum_{i=1}^{k_n}\xi_{n,i},\quad\bar{\xi}_n^*:=k_n^{-1}\sum_{i=1}^{k_n}\xi_{n,i}^*, \quad\mbox{and}\quad\varsigma_n^2:=\Var_{\prob}\xi_{n,1}.
\]
If
\begin{equation}\label{CLTbootIID2}
\liminf_{n\to\infty}\varsigma_n^2=\varsigma^2>0,
\end{equation}
then
\begin{equation*}
\sup_{x\in\mathbbm{R}}\left|\prob_{\bxi}^*\left[\frac{\sqrt{k_n}}{\sqrt{\varsigma_n^2}}\left(\bar{\xi}_n^*-\bar{\xi}_n\right)\leq x\right]-\prob\left[\frac{\sqrt{k_n}}{\sqrt{\varsigma_n^2}}\bar{\xi}_n\leq x\right]\right|\xrightarrow[n\to\infty]{\prob}0.
\end{equation*}
\end{theorem}

\begin{theorem}[Bootstrap multivariate CLT for triangular arrays]\label{thm:CLTbootTRIAmult}
Let $\{\bxi_{n,k_n}\}_{n=1}^{\infty}$ be a~triangular array of zero mean $q$-dimensional random vectors on the same probability space such that the elements of the vector sequence $\{\bxi_{n,1},\ldots,\bxi_{n,k_n}\}$ are iid for every $n\in\mathbbm{N}$ satisfying
\begin{equation}\label{CLTbootIID1mult}
\sup_{n\in\mathbbm{N}}\E_{\prob}|\xi_{n,1}^{(j)}|^4<\infty,\quad j\in\{1,\ldots,q\},
\end{equation}
where $\bxi_{n,1}\equiv[\xi_{n,1}^{(1)},\ldots,\xi_{n,1}^{(q)}]^{\top}\in\mathbbm{R}^q,\,n\in\mathbbm{N}$ and $k_n\to\infty$ as $n\to\infty$. Assume that $\bXi^*\equiv[\bxi_{n,1}^*,\ldots,\bxi_{n,k_n}^*]^{\top}$ is the bootstrapped version of $\bXi\equiv[\bxi_{n,1},\ldots,\bxi_{n,k_n}]^{\top}$. Denote
\[
\bar{\bxi}_n:=k_n^{-1}\sum_{i=1}^{k_n}\bxi_{n,i},\quad\bar{\bxi}_n^*:=k_n^{-1}\sum_{i=1}^{k_n}\bxi_{n,i}^*, \quad\mbox{and}\quad\bGamma_n:=\Var_{\prob}\bxi_{n,1}.
\]
If
\begin{equation}\label{CLTbootIID2mult}
\liminf_{n\to\infty}\bGamma_n=\bGamma>\zero,
\end{equation}
then
\begin{equation*}
\prob_{\bXi}^*\left[\sqrt{k_n}\bGamma_n^{-1/2}\left(\bar{\bxi}_n^*-\bar{\bxi}_n\right)\leq \bx\right]-\prob\left[\sqrt{k_n}\bGamma_n^{-1/2}\bar{\bxi}_n\leq \bx\right]\xrightarrow[n\to\infty]{\prob}0,\quad\forall\bx\in\mathbbm{R}^q.
\end{equation*}
\end{theorem}

\section{Proofs}

\begin{proof}[of Theorem~\ref{underNull}]
Let us define
$$
U_N(t):=\frac{1}{\sigma\sqrt{N}}\sum_{i=1}^N\sum_{s=1}^t(Y_{i,s}-\mu_i).
$$
Using the multivariate Lindeberg-L\'{e}vy CLT for a~sequence of $T$-dimensional iid random vectors $\{[\sum_{s=1}^1\eps_{i,s},\ldots,\sum_{s=1}^T\eps_{i,s}]^{\top}\}_{i\in\mathbbm{N}}$, we have under $H_0$
$$
[U_N(1),\ldots,U_N(T)]^{\top}\xrightarrow[N\to\infty]{\dist}[X_1,\ldots,X_T]^{\top},
$$
since $\Var[\sum_{s=1}^1\eps_{1,s},\ldots,\sum_{s=1}^T\eps_{1,s}]^{\top}=\bLambda$. Indeed, the $t$-th diagonal element of the covariance matrix $\bLambda$ is
$$
\Var\sum_{s=1}^t\eps_{1,s}=r(t)
$$
and the upper off-diagonal element on position $(t,v)$ is
\begin{align*}
\Cov\left(\sum_{s=1}^t\eps_{1,s},\sum_{u=1}^v\eps_{1,u}\right)&=\Var\sum_{s=1}^t\eps_{1,s}+\Cov\left(\sum_{s=1}^t\eps_{1,s},\sum_{u=t+1}^v\eps_{1,u}\right)\\
&=r(t)+R(t,v),\quad t<v.
\end{align*}
Moreover, let us define the reverse analogue to $U_N(t)$, i.e.,
$$
V_N(t):=\frac{1}{\sigma\sqrt{N}}\sum_{i=1}^N\sum_{s=t+1}^T(Y_{i,s}-\mu_i)=U_N(T)-U_N(t).
$$
Hence,
\begin{align*}
U_{N}(s)-\frac{s}{t}U_{N}(t)&=\frac{1}{\sigma\sqrt{N}}\sum_{i=1}^N\left\{\sum_{r=1}^s\left[\left(Y_{i,r}-\mu_i\right)-\frac{1}{t}\sum_{v=1}^t \left(Y_{i,v}-\mu_i\right)\right]\right\}\\
&=\frac{1}{\sigma\sqrt{N}}\sum_{i=1}^N\sum_{r=1}^s\left(Y_{i,r}-\widebar{Y}_{i,t}\right)
\end{align*}
and, consequently,
\begin{align*}
V_{N}(s)-\frac{T-s}{T-t}V_{N}(t)&=\frac{1}{\sigma\sqrt{N}}\sum_{i=1}^N\left\{\sum_{r=s+1}^T\left[\left(Y_{i,r}-\mu_i\right)-\frac{1}{T-t}\sum_{v=t+1}^T \left(Y_{i,v}-\mu_i\right)\right]\right\}\\
&=\frac{1}{\sigma\sqrt{N}}\sum_{i=1}^N\sum_{r=s+1}^T\left(Y_{i,r}-\widetilde{Y}_{i,t}\right).
\end{align*}
Using the Cram\'{e}r-Wold device, we end up with
\begin{multline*}
\max_{t=2,\ldots,T-2}\frac{\max_{s=1,\ldots,t}\left|U_{N}(s)-\frac{s}{t}U_{N}(t)\right|}{\max_{s=t,\ldots,T-1}\left|V_{N}(s)-\frac{T-s}{T-t}V_{N}(t)\right|}\\
\xrightarrow[N\to\infty]{\dist}\max_{t=2,\ldots,T-2}\frac{\max_{s=1,\ldots,t}\left|X_s-\frac{s}{t}X_t\right|}{\max_{s=t,\ldots,T-1}\left|(X_T-X_s)-\frac{T-s}{T-t}(X_T-X_t)\right|}.
\end{multline*}
\qed
\end{proof}

\begin{proof}[of Theorem~\ref{underAlternative}]
Let $t=\tau+1$. Then, under alternative $H_1$
\begin{align*}
&\frac{1}{\sigma\sqrt{N}}\max_{s=1,\ldots,\tau+1}\left|\sum_{i=1}^N\left[\sum_{r=1}^s\left(Y_{i,r}-\widebar{Y}_{i,\tau+1}\right)\right]\right|\nonumber\\
&\geq\frac{1}{\sigma\sqrt{N}}\left|\sum_{i=1}^N\sum_{r=1}^{\tau}\left(Y_{i,r}-\widebar{Y}_{i,\tau+1}\right)\right|\\
&=\frac{1}{\sigma\sqrt{N}}\left|\sum_{i=1}^N\sum_{r=1}^{\tau}\left(\mu_i+\sigma\eps_{i,r}-\frac{1}{\tau+1}\sum_{v=1}^{\tau+1}(\mu_i+\sigma\eps_{i,v})-\frac{1}{\tau+1}\delta_i\right)\right|\\
&=\frac{1}{\sqrt{N}}\left|\sum_{i=1}^N\sum_{r=1}^{\tau}\left(\eps_{i,r}-\widebar{\eps}_{i,\tau+1}\right)-\frac{\tau}{\sigma(\tau+1)}\sum_{i=1}^N\delta_i\right|\\
&=\mathcal{O}_{\prob}(1)+\frac{\tau}{\sigma(\tau+1)\sqrt{N}}\left|\sum_{i=1}^N\delta_i\right|\xrightarrow[]{\prob}\infty,\quad N\to\infty,
\end{align*}
where $\widebar{\eps}_{i,\tau+1}=\frac{1}{\tau}\sum_{v=1}^{\tau+1}\eps_{i,v}$.

Since there is no change after $\tau+1$ and $\tau\leq T-3$, then by Theorem~\ref{underNull} we have
$$
\frac{1}{\sigma\sqrt{N}}\max_{s=\tau+1,\ldots,T-1}\left|\sum_{i=1}^N\sum_{r=s+1}^T\left(Y_{i,r}-\widetilde{Y}_{i,\tau+1}\right)\right|\xrightarrow[N\to\infty]{\dist}\max_{s=\tau+1,\ldots,T-1}\left|Z_s-\frac{T-s}{T-\tau}Z_{\tau+1}\right|.
$$
\qed
\end{proof}

\begin{proof}[of Theorem~\ref{thm:change}]
Let us define
\[
S_N^{(i)}(t):=\frac{1}{w(t)}\sum_{s=1}^t(Y_{i,s}-\widebar{Y}_{i,t})^2
\]
and, consequently, $S_N(t):=\frac{1}{N}\sum_{i=1}^N S_N^{(i)}(t)$. Then,
\[
S_N^{(i)}(t)=\left\{\begin{array}{ll}
\frac{\sigma^2}{w(t)}\sum_{s=1}^t(\e_{i,s}-\widebar{\e}_{i,t})^2, & t\leq\tau,\\
\frac{1}{w(t)}\Big[\sum_{s=1}^{\tau}(\sigma\e_{i,s}-\sigma\widebar{\e}_{i,t}-\frac{t-\tau}{t}\delta_i)^2\Big.&\\
\quad\Big.+\sum_{s=\tau+1}^{t}(\sigma\e_{i,s}-\sigma\widebar{\e}_{i,t}+\frac{\tau}{t}\delta_i)^2\Big], & t>\tau;
\end{array}
\right.
\]
where $\widebar{\e}_{i,t}=\frac{1}{t}\sum_{s=1}^t\e_{i,s}$. By the definition of the cumulative autocorrelation function, we have for $2\leq t\leq\tau$
\begin{align*}
\E S_N^{(i)}(t)&=\frac{\sigma^2}{w(t)}\sum_{s=1}^t \E(\e_{i,s}-\widebar{\e}_{i,t})^2=\frac{\sigma^2}{w(t)}\sum_{s=1}^t\left[1-\frac{2}{t}\sum_{r=1}^t\E\e_{i,s}\e_{i,r}+\frac{1}{t^2}r(t)\right]\\
&=\frac{\sigma^2}{w(t)}\left(t-\frac{r(t)}{t}\right).
\end{align*}
In the other case when $t>\tau$, one can calculate
\begin{align*}
\E S_N^{(i)}(t)&=\frac{\sigma^2}{w(t)}\left(t-\frac{r(t)}{t}\right)+\frac{\tau}{w(t)}\left(\frac{t-\tau}{t}\right)^2\delta_i^2+\frac{t-\tau}{w(t)}\left(\frac{\tau}{t}\right)^2\delta_i^2\\
&=\frac{\sigma^2t}{w(t)}\left(1-\frac{r(t)}{t^2}\right)+\frac{\tau(t-\tau)}{tw(t)}\delta_i^2.
\end{align*}
Realize that $S_N^{(i)}(t)-\E S_N^{(i)}(t)$ are independent with zero mean for fixed $t$ and $i=1,\ldots,N$. Due to Assumption~\ref{ass:four}, for $2\leq t\leq\tau$ it holds
\[
\Var S_N(t)=\frac{1}{N^2}\sum_{i=1}^N\frac{\sigma^4}{w^2(t)}\Var \left[\sum_{s=1}^t(\e_{i,s}-\widebar{\e}_{i,t})^2\right]=\frac{1}{N}C_1(t,\sigma),
\]
where $C_1(t,\sigma)>0$ is some constant not depending on $N$. If $t>\tau$, then
\begin{align*}
\Var S_N(t)&=\frac{1}{N^2}\sum_{i=1}^N\frac{1}{w^2(t)}\Var\Bigg[\sigma^2\sum_{s=1}^{\tau}(\e_{i,s}-\widebar{\e}_{i,t})^2-2\frac{t-\tau}{t}\sigma\delta_i\sum_{s=1}^{\tau}(\e_{i,s}-\widebar{\e}_{i,t})\Bigg.\\
&\quad+\Bigg.\left(\frac{t-\tau}{t}\right)^2\delta_i^2+\sigma^2\sum_{s=\tau+1}^t(\e_{i,s}-\widebar{\e}_{i,t})^2\\
&\quad+2\frac{\tau}{t}\sigma\delta_i\sum_{s=\tau+1}^t(\e_{i,s}-\widebar{\e}_{i,t})+\left(\frac{\tau}{t}\right)^2\delta_i^2\Bigg]\\
&\leq \frac{1}{N}C_2(t,\tau,\sigma)+\frac{1}{N^2}C_3(t,\tau,\sigma)\sum_{i=1}^N\delta_i^2+\frac{1}{N^2}C_4(t,\tau,\sigma)\left|\sum_{i=1}^N\delta_i\right|,
\end{align*}
where $C_j(t,\tau,\sigma)>0$ does not depend on $N$ for $j=2,3,4$.

The Chebyshev inequality provides $S_N(t)-\E S_N(t)=\Op\left(\sqrt{\Var S_N(t)}\right)$ as $N\to\infty$. According to Assumption~\ref{ass:tozero} and the Cauchy-Schwarz inequality, we have
\[
\frac{1}{N^2}\left|\sum_{i=1}^N\delta_i\right|\leq \frac{1}{N}\sqrt{\frac{1}{N}\sum_{i=1}^N\delta_i^2}\to 0,\quad N\to\infty.
\]
Since the index set $\{2,\ldots,T\}$ is finite and $\tau$ is finite as well, then
\[
\max_{2\leq t\leq T}\Var S_N(t)\leq \frac{1}{N}K_1(\sigma)+K_2(\sigma)\frac{1}{N^2}\sum_{i=1}^N\delta_i^2+K_3(\sigma)\frac{1}{N^2}\left|\sum_{i=1}^N\delta_i\right|\leq\frac{1}{N}K_4(\sigma),
\]
where $K_j(\sigma)>0$ are constants not depending on $N$ for $j=1,2,3,4$. Thus, we also have uniform stochastic boundedness, i.e.,
\[
\max_{2\leq t\leq T}|S_N(t)-\E S_N(t)|=\Op\left(\frac{1}{\sqrt{N}}\right),\quad N\to\infty.
\]

Adding and subtracting, one has
\begin{align*}
S_N(t)-S_N(\tau)&=S_N(t)-\E S_N(t)-[S_N(\tau)-\E S_N(\tau)]\\
&\quad+\E S_N(t)-\E S_N(\tau)\\
&\geq -2\max_{2\leq r\leq T}|S_N(r)-\E S_N(r)|+\E S_N(t)-\E S_N(\tau)\\
&= -2\max_{2\leq r\leq T}|S_N(r)-\E S_N(r)|\\
&\quad+\sigma^2\left[\frac{t}{w(t)}\left(1-\frac{r(t)}{t^2}\right)-\frac{\tau}{w(\tau)}\left(1-\frac{r(\tau)}{\tau^2}\right)\right]\\
&\quad+\mathcal{I}\{t>\tau\}\frac{\tau(t-\tau)}{tw(t)}\frac{1}{N}\sum_{i=1}^N\delta_i^2.
\end{align*}
The above inequality holds for each $t\in\{2,\ldots,T\}$ and, particularly, it holds for $\widehat{\tau}_N$. Note that $\widehat{\tau}_N=\arg\min_tS_N(t)$. Hence, $S_N(\widehat{\tau}_N)-S_N(\tau)\leq 0$. Therefore,
\begin{align}\label{changeineq}
&2\sqrt{N}\max_{2\leq r\leq T}|S_N(r)-\E S_N(r)|\nonumber\\
&\quad\geq \sqrt{N}\Bigg\{\sigma^2\left[\frac{\widehat{\tau}_N}{w(\widehat{\tau}_N)}\left(1-\frac{r(\widehat{\tau}_N)}{\widehat{\tau}_N^2}\right)-\frac{\tau}{w(\tau)}\left(1-\frac{r(\tau)}{\tau^2}\right)\right]\Bigg.\nonumber\\
&\quad+\Bigg.\mathcal{I}\{\widehat{\tau}_N>\tau\}\frac{\tau(\widehat{\tau}_N-\tau)}{\widehat{\tau}_Nw(\widehat{\tau}_N)}\frac{1}{N}\sum_{i=1}^N\delta_i^2\Bigg\}.
\end{align}
If $\widehat{\tau}_N>\tau$ almost surely for infinitely many $N$, then the left hand side of~\eqref{changeineq} is $\Op(1)$ as $N\to\infty$, but the right hand side is unbounded because of Assumption~\ref{ass:delta2E}. If $\widehat{\tau}_N<\tau$ almost surely for infinitely many $N$, then due to the monotonicity Assumption~\ref{ass:decreasing}
\begin{multline*}
0\xleftarrow[N\to\infty]{\P} 2\max_{2\leq r\leq T}|S_N(r)-\E S_N(r)|\\
\geq \sigma^2\left[\frac{\widehat{\tau}_N}{w(\widehat{\tau}_N)}\left(1-\frac{r(\widehat{\tau}_N)}{\widehat{\tau}_N^2}\right)-\frac{\tau}{w(\tau)}\left(1-\frac{r(\tau)}{\tau^2}\right)\right]>0,
\end{multline*}
which is a~contradicting conclusion. Hence, $\P[\widehat{\tau}_N=\tau]\to 1$ as $N\to\infty$.
\end{proof}

\begin{proof}[of Theorem~\ref{bootJust}]
Let us define $\widehat{\epsilon}_{i,t}:=\sigma^{-1}\sum_{s=1}^t\widehat{e}_{i,s}$, $\widehat{\epsilon}_{i,t}^*:=\sigma^{-1}\sum_{s=1}^t\widehat{e}_{i,s}^*$,
\[
\widehat{U}_N(t):=\frac{1}{\sigma\sqrt{N}}\sum_{i=1}^N\sum_{s=1}^t\widehat{e}_{i,s}=\frac{1}{\sqrt{N}}\sum_{i=1}^N\widehat{\epsilon}_{i,t},
\]
and
\begin{align*}
\widehat{U}_N^*(t)&:=\frac{1}{\sigma\sqrt{N}}\sum_{i=1}^N\sum_{s=1}^t\widehat{Y}_{i,s}^*=\frac{1}{\sigma\sqrt{N}}\sum_{i=1}^N\sum_{s=1}^t\left(\widehat{e}_{i,s}^*-\frac{1}{N}\sum_{i=1}^N\widehat{e}_{i,s}\right)\\
&=\frac{1}{\sigma\sqrt{N}}\sum_{i=1}^N\sum_{s=1}^t\left(\widehat{e}_{i,s}^*-\widehat{e}_{i,s}\right)=\frac{1}{\sqrt{N}}\sum_{i=1}^N\left(\widehat{\epsilon}_{i,t}^*-\widehat{\epsilon}_{i,t}\right).
\end{align*}

Realize that $\widehat{\epsilon}_{i,t}$ depends on $\widehat{\tau}_N$ and, hence, it depends on $N$. Thus, $\widehat{\epsilon}_{i,t}\equiv\widehat{\epsilon}_{i,t}(N)$. Since Assumption~\ref{ass:four} holds, then according to the bootstrap multivariate CLT for triangular arrays (Theorem~\ref{thm:CLTbootTRIAmult}) of $T$-dimensional vectors $\bxi_{N,i}=[\widehat{\epsilon}_{i,1}(N),\ldots,\widehat{\epsilon}_{i,T}(N)]^{\top}$ with $k_N=N$, we have
\begin{multline*}
\prob\left[\bGamma_N^{-1/2}[\widehat{U}_N^*(1),\ldots,\widehat{U}_N^*(T)]^{\top}\leq \bx\big|\mathbbm{Y}\right]-\prob\left[\bGamma_N^{-1/2}[\widehat{U}_N(1),\ldots,\widehat{U}_N(T)]^{\top}\leq \bx\right]\\
\xrightarrow[N\to\infty]{\prob}0,\quad\forall\bx\in\mathbbm{R}^T,
\end{multline*}
where $\bGamma_N=\Var[\widehat{\epsilon}_{i,1},\ldots,\widehat{\epsilon}_{i,T}]^{\top}$.

Now, it is sufficient to realize that $[\widehat{U}_N(1),\ldots,\widehat{U}_N(T)]^{\top}$ has an~approximate multivariate normal distribution with zero mean and covariance matrix $\bGamma=\lim_{N\to\infty}\bGamma_N$. Using the law of total variance,
\begin{align*}
\Var\widehat{\epsilon}_{i,t}&=\E[\Var\{\widehat{\epsilon}_{i,t}|\widehat{\tau}_N\}]+\Var[\E\{\widehat{\epsilon}_{i,t}|\widehat{\tau}_N\}]\\
&=\sum_{\pi=1}^T\prob[\widehat{\tau}_N=\pi]\Var[\widehat{\epsilon}_{i,t}|\widehat{\tau}_N=\pi]+\sum_{\pi=1}^T\prob[\widehat{\tau}_N=\pi]\{\E[\widehat{\epsilon}_{i,t}|\widehat{\tau}_N=\pi]\}^2\\
&\quad-\left\{\sum_{\pi=1}^T\prob[\widehat{\tau}_N=\pi]\E[\widehat{\epsilon}_{i,t}|\widehat{\tau}_N=\pi]\right\}^2.
\end{align*}
Since $\lim_{N\to\infty}\prob[\widehat{\tau}_N=\tau]=1$ and $\E[\widehat{e}_{i,t}|\widehat{\tau}_N=\tau]=0$, then
\[
\lim_{N\to\infty}\Var\widehat{\epsilon}_{i,t}=\lim_{N\to\infty}\Var[\widehat{\epsilon}_{i,t}|\widehat{\tau}_N=\tau].
\]
Similarly with the covariance, i.e., after applying the law of total covariance, we have
\[
\lim_{N\to\infty}\Cov\left(\widehat{\epsilon}_{i,t},\widehat{\epsilon}_{i,v}\right)=\lim_{N\to\infty}\Cov\left(\widehat{\epsilon}_{i,t},\widehat{\epsilon}_{i,v}|\widehat{\tau}_N=\tau\right).
\]
Note that
\[
\left(\widehat{e}_{i,t}|\widehat{\tau}_N=\tau\right)=\left\{
\begin{array}{ll}
\sigma(\eps_{i,t}-\widebar{\eps}_{i,\tau}),& t\leq\tau;\\
\sigma(\eps_{i,t}-\widetilde{\eps}_{i,\tau}),& t>\tau;
\end{array}
\right.
\]
where
\[
\widebar{\eps}_{i,t}=\frac{1}{t}\sum_{s=1}^t \eps_{i,s}\quad\mbox{and}\quad\widetilde{\eps}_{i,t}=\frac{1}{T-t}\sum_{s=t+1}^T \eps_{i,s}.
\]
Taking into account the definitions of $r(t)$, $R(t,v)$, and $S(t,v,d)$ together with some simple algebra, we obtain that $\Var[\widehat{\epsilon}_{i,s}|\widehat{\tau}_N=\tau]=\gamma_{t,t}(\tau)$ and $\Cov\left(\widehat{\epsilon}_{i,t},\widehat{\epsilon}_{i,v}|\widehat{\tau}_N=\tau\right)=\gamma_{t,v}(\tau)$ for $t<v$, where the elements $\gamma_{t,t}(\tau)$ and $\gamma_{t,v}(\tau)$ are as in the statement of Theorem~\ref{bootJust}.

Then the sum in the nominator of $\mathcal{R}_N^*(T)$ can be alternatively rewritten as
\[
\frac{1}{\sigma\sqrt{N}}\sum_{i=1}^N\sum_{r=1}^s\left(\widehat{Y}_{i,r}^*-\widebar{\widehat{Y}}_{i,t}^*\right)=\frac{1}{\sigma\sqrt{N}}\sum_{i=1}^N\left\{\left[\sum_{r=1}^s \widehat{Y}_{i,r}^*\right]-\frac{s}{t}\sum_{v=1}^t \widehat{Y}_{i,v}^*\right\}=\widehat{U}_{N}^*(s)-\frac{s}{t}\widehat{U}_{N}^*(t).
\]

Concerning the denominator of $\mathcal{R}_N^*(T)$, one needs to perform a~similar calculation as in the proof of Theorem~\ref{underNull} with $V_N(t)$, i.e., to define $\widehat{V}_N(t)$ and $\widehat{V}_N^*(t)$ analogously to $\widehat{U}_N(t)$ and $\widehat{U}_N^*(t)$ as $V_N(t)$ is to $U_N(t)$. Applying the Cram\'{e}r-Wold theorem completes the proof.
\qed
\end{proof}

\begin{proof}[of Theorem~\ref{bootConsis}]
Recall the notation from the proof of Theorem~\ref{bootJust}. Under $H_0$, \ref{withprob}, and~\ref{ass:four} it holds
\[
\lim_{N\to\infty}\prob[\widehat{\tau}_N=T]=1.
\]
Then in view of~\eqref{widehate},
\[
\lim_{N\to\infty}\prob\left[\widehat{U}_N(s)-\frac{s}{t}\widehat{U}_N(t)=U_N(s)-\frac{s}{t}U_N(t)\right]=1,\quad 1\leq s\leq t\leq T.
\]
\qed
\end{proof}


\begin{proof}[of Theorem~\ref{thm:CLTbootTRIA}]
The Lyapunov condition~\cite[p.~371]{Billingsley1986} for a~triangular array of random variables $\{\xi_{n,k_n}\}_{n=1}^{\infty}$ is satisfied due to~\eqref{CLTbootIID1} and~\eqref{CLTbootIID2}, i.e., for $\omega=2$:
\[
\frac{1}{\sqrt{k_n\varsigma_n^2}^{2+\omega}}\sum_{i=1}^{k_n}\E|\xi_{n,i}|^{2+\omega}\leq\frac{k_n^{-\omega/2}}{\varsigma_n^{2+\omega}}\sup_{\iota\in\mathbbm{N}}\E|\xi_{\iota,1}|^{2+\omega}\to 0,\quad n\to\infty.
\]
Therefor, the CLT for $\{\xi_{n,k_n}\}_{n=1}^{\infty}$ holds and
\[
\sup_{x\in\mathbbm{R}}\left|\prob\left[\frac{\sqrt{k_n}}{\sqrt{\varsigma_n^2}}\bar{\xi}_n\leq x\right]-\int_{-\infty}^x\frac{1}{\sqrt{2\pi}}\exp\left\{-\frac{t^2}{2}\right\}\ud t\right|\xrightarrow[n\to\infty]{}0.
\]
Now, to prove the theorem, it suffices to show the following three statements:
\begin{enumerate}\setlength{\itemsep}{1mm}
\item $\sup_{x\in\mathbbm{R}}\left|\prob_{\bxi}^*\left[\frac{\sqrt{k_n}}{\sqrt{\Var_{\prob_{\bxi}^*}\xi_{n,1}^*}}\left(\bar{\xi}_n^*-\E_{\prob_{\bxi}^*}\bar{\xi}_n^*\right)\leq x\right]-\int_{-\infty}^x\frac{1}{\sqrt{2\pi}}\exp\left\{-\frac{t^2}{2}\right\}\ud t\right|\xrightarrow[n\to\infty]{\prob}0$;
\item $\Var_{\prob_{\bxi}^*}\xi_{n,1}^*-\varsigma_n^2\xrightarrow[n\to\infty]{\prob}0$;
\item $\E_{\prob_{\bxi}^*}\bar{\xi}_n^*=\bar{\xi}_n,\, [\prob]\mbox{-}a.s.$
\end{enumerate}

Proving~\emph{(iii)} is trivial, because $\E_{\prob_{\bxi}^*}\bar{\xi}_n^*=\E_{\prob_{\bxi}^*}\xi_{n,1}^*=k_n^{-1}\sum_{i=1}^{k_n}\xi_{n,i}=\bar{\xi}_n,\, [\prob]\mbox{-}a.s.$

Let us calculate the conditional variance of the bootstrapped variable $\xi_{n,1}^*$: $\Var_{\prob_{\bxi}^*}\xi_{n,1}^*=\E_{\prob_{\bxi}^*}\xi_{n,1}^{*2}-(\E_{\prob_{\bxi}^*}\xi_{n,1}^*)^2=k_n^{-1}\sum_{i=1}^{k_n}\xi_{n,i}^2-\left(k_n^{-1}\sum_{i=1}^{k_n}\xi_{n,i}\right)^2,\, [\prob]\mbox{-}a.s.$ The weak law of large numbers together with~\eqref{CLTbootIID1} provides
\[
\bar{\xi}_n-n^{-1}\sum_{i=i}^{k_n}\E_{\prob}\xi_{n,i}=\bar{\xi}_n\xrightarrow[n\to\infty]{\prob}0
\]
and
\[
0\xleftarrow[n\to\infty]{\prob}k_n^{-1}\sum_{i=1}^{k_n}\xi_{n,i}^2-\left(k_n^{-1}\sum_{i=1}^{k_n}\xi_{n,i}\right)^2-\E_{\prob}\xi_{n,1}^2=\Var_{\prob_{\bxi}^*}\xi_{n,1}^*-\varsigma_n^2.
\]
The last result of the WLLN is true, because~\eqref{CLTbootIID1} implies
\[
k_n^{-2}\sum_{i=1}^{k_n}\Var_{\prob}\xi_{n,i}^2\leq k_n^{-2}\sum_{i=1}^{k_n}\E_{\prob}\xi_{n,i}^4\leq k_n^{-1}\sup_{\iota\in\mathbbm{N}}\E_{\prob}\xi_{\iota,1}^4\xrightarrow[n\to\infty]{} 0.
\]
Thus \emph{(ii)} is proved.

The Berry-Esseen-Katz theorem (see~\cite{Katz1963}) with $g(x)=|x|^{\epsilon},\,\epsilon>0$ for the bootstrapped sequence of $iid$ (with respect to $\prob^*$) random variables $\{\xi_{n,i}^*\}_{i=1}^{k_n}$ results in
\begin{multline}\label{bekatz}
\sup_{x\in\mathbbm{R}}\left|\prob_{\bxi}^*\left[\frac{\sqrt{k_n}}{\sqrt{\Var_{\prob_{\bxi}^*}\xi_{n,1}^*}}\left(\bar{\xi}_n^*-\E_{\prob_{\bxi}^*}\bar{\xi}_n^*\right)\leq x\right]-\int_{-\infty}^x\frac{1}{\sqrt{2\pi}}\exp\left\{-\frac{t^2}{2}\right\}\ud t\right|\\
\leq Ck_n^{-\epsilon/2}\E_{\prob_{\bxi}^*}\left|\frac{\xi_{n,1}^*-\E_{\prob_{\bxi}^*}\xi_{n,1}^*}{\Var_{\prob_{\bxi}^*}\xi_{n,1}^*}\right|^{2+\epsilon}\quad [\prob]\mbox{-}a.s.,
\end{multline}
for all $n\in\mathbbm{N}$ where $C>0$ is an~absolute constant.

The Jensen inequality and Minkowski inequality provide an~upper bound for the nominator from the right-hand side of~\eqref{bekatz}:
\begin{align*}
&\E_{\prob_{\bxi}^*}|\xi_{n,1}^*-\E_{\prob_{\bxi}^*}\xi_{n,1}^*|^{2+\epsilon}=k_n^{-1}\sum_{i=1}^{k_n}\left|\xi_{n,i}-k_n^{-1}\sum_{j=1}^{k_n}\xi_{n,j}\right|^{2+\epsilon}\\
&\leq k_n^{-1}\left\{\left(\sum_{i=1}^{k_n}|\xi_{n,i}|^{2+\epsilon}\right)^{1/(2+\epsilon)}+k_n^{-(1+\epsilon)/(2+\epsilon)}\left|\sum_{j=1}^{k_n}\xi_{n,j}\right|\right\}^{2+\epsilon}\\
&\leq 2^{1+\epsilon}k_n^{-1}\sum_{i=1}^{k_n}|\xi_{n,i}|^{2+\epsilon}+2^{1+\epsilon}\left| k_n^{-1}\sum_{i=1}^{k_n}\xi_{n,i}\right|^{2+\epsilon}\quad [\prob]\mbox{-}a.s.
\end{align*}
The right-hand side of the previously derived upper bound is uniformly bounded in probability $\prob$, because of Markov's inequality and~\eqref{CLTbootIID1}. Indeed, for fixed $\eta>0$
\[
\prob\left[k_n^{-1}\sum_{i=1}^{k_n}|\xi_{n,i}|^{2+\epsilon}\geq\eta\right]\leq \eta^{-1}k_n^{-1}\sum_{i=1}^{k_n}\E_{\prob}|\xi_{n,i}|^{2+\epsilon}\leq\eta^{-1}\sup_{\iota\in\mathbbm{N}}\E_{\prob}|\xi_{\iota,1}|^{2+\epsilon}<\infty,\quad\forall n\in\mathbbm{N}
\]
and
\[
\prob\left[\left| k_n^{-1}\sum_{i=1}^{k_n}\xi_{n,i}\right|\geq\eta\right]\leq \eta^{-1}k_n^{-1}\E_{\prob}\left|\sum_{i=1}^{k_n}\xi_{n,i}\right|\leq\eta^{-1}\sup_{\iota\in\mathbbm{N}}\E_{\prob}|\xi_{\iota,1}|<\infty,\quad\forall n\in\mathbbm{N}.
\]
Since $\E_{\prob_{\bxi}^*}|\xi_{n,1}^*-\E_{\prob^*}\xi_{n,1}^*|^{2+\epsilon}$ is bounded in probability $\prob$ uniformly over $n$ and the denominator of the right-hand side of~\eqref{bekatz} is uniformly bounded away from zero due to~\eqref{CLTbootIID2}, then the left-hand side of~\eqref{bekatz} converges in probability $\prob$ to zero as $n$ tends to infinity. So, \emph{(i)} is proved as well.
\qed
\end{proof}

\begin{proof}[of Theorem~\ref{thm:CLTbootTRIAmult}]
According to the Cram\'{e}r-Wold theorem, it is sufficient to ensure that all assumptions of one-dimensional bootstrap CLT~\ref{thm:CLTbootTRIA} for triangular arrays are valid for any linear combination of the elements of the random vector $\bxi_{n,1},\,n\in\mathbbm{N}$.

For arbitrary fixed $\bt\in\mathbbm{R}^q$ using the Jensen inequality, we get
\[
\sup_{n\in\mathbbm{N}}\E_{\prob}|\bt^{\top}\bxi_{n,1}|^4\leq q^3\sup_{n\in\mathbbm{N}}\sum_{j=1}^q t_j^4\E_{\prob}|\xi_{n,1}^{(j)}|^4\leq q^4\max_{j=1,\ldots,q}t_j^4\sup_{n\in\mathbbm{N}}\E_{\prob}|\xi_{n,1}^{(j)}|^4<\infty.
\]
Hence, assumption~\eqref{CLTbootIID1mult} implies assumption~\eqref{CLTbootIID1} for the random variables $\{\bt^{\top}\bxi_{n,k_n}\}_{n\in\mathbbm{N}}$.

Similarly, assumption~\eqref{CLTbootIID2mult} implies assumption~\eqref{CLTbootIID2} for such an~arbitrary linear combination, i.e., positive definiteness of the matrix $\bGamma$ yields
\[
\liminf_{n\to\infty}\Var_{\prob}\bt^{\top}\bxi_{n,1}=\liminf_{n\to\infty}\bt^{\top}\left(\Var_{\prob}\bxi_{n,1}\right)\bt\geq \bt^{\top}\left(\liminf_{n\to\infty}\bGamma_n\right)\bt=\bt^{\top}\bGamma\bt>0.
\]
\qed
\end{proof}

\begin{acknowledgements}
The authors thank two anonymous referees and the Associate Editor for the suggestions that improved the paper. The authors also thank Professor Daniela Jaru\v{s}kov\'{a} and Professor Zuzana Pr\'{a}\v{s}kov\'{a} for pointing out the mistake in the original paper. This corrected paper was written with the support of the Czech Science Foundation project GA\v{C}R No.~P201/13/12994P.
\end{acknowledgements}

\bibliographystyle{spmpsci}      
\bibliography{pesta}   


\end{document}